\documentclass{article}
\usepackage{enumerate,cancel}
\usepackage{amssymb,amsfonts,amsmath,latexsym,color, amscd} 
\usepackage[dvips]{psfrag,graphicx} 
\usepackage[latin1]{inputenc}

\def\RR{\mathbb R}

\newcommand{\EE}{\mathbb{E}}

\def\ss{{\mathbb{S}}}

\def\si{\sigma}
\def\ga{\gamma}

 \def\ppt#1{{\overset{\mbox{\large .\kern-1pt.}}{#1}}}
 \def\pppt#1{{\overset{\mbox{\large .\kern-1pt.\kern-1pt.}}{#1}}}

 \def\Pt#1{{\overset{\mbox{\Huge .}}{#1}}}
 \def\Ppt#1{{\overset{\mbox{\Huge .\kern-3pt.}}{#1}}}
 \def\Pppt#1{{\overset{\mbox{\Huge .\kern-3pt.\kern-3pt.}}{#1}}}

\def\c#1{\overset{\mbox{\tiny $\circ$}}{#1}}
\def\cc#1{\overset{\mbox{\tiny $\circ\circ$}}{#1}}
\def\ccc#1{\overset{\mbox{\tiny $\circ\!\!\circ\!\!\circ$}}{#1}}

\def\vector#1{\mbox{\boldmath $#1$}}
\def\vect#1{{\overset{\rightarrow}{#1}}}

\renewcommand{\phi}{\varphi}
\renewcommand{\epsilon}{\varepsilon}

\newcommand{\esc}[2]{\langle #1,#2\rangle}

\renewcommand{\epsilon}{\varepsilon}

\newcommand{\spa}{\mathrm{span}}

\newtheorem{theo}{Theorem}[section]
\newtheorem{lemm}[theo]{Lemma}
\newtheorem{coro}[theo]{Corollary}
\newtheorem{defi}[theo]{Definition}
\newtheorem{prop}[theo]{Proposition}

\newtheorem{remark}[theo]{Remark}

\newenvironment{demo}{\noindent {\bf Proof: }}{\hfill$\Box$\par\bigskip}

\numberwithin{equation}{section}

\title{Space of subspheres and conformal invariants of curves}
\author{R\'{e}mi Langevin, \footnote{Institut de Math\'ematiques de Bourgogne, UMR C.N.R.S. 5584, Universit\'e de Bourgogne}\\Jun O'Hara and Shigehiro Sakata \footnote{Department of Mathematics, Tokyo Metropolitan University}}

\begin{document}\thispagestyle{empty}

\maketitle

\begin{abstract}{A space curve is determined by conformal arc-length, conformal curvature, and conformal torsion, up to M\"obius transformations. 
We use the spaces of osculating circles and spheres to give a conformally defined moving frame of a curve in the Minkowski space, which can naturally produce the conformal invariants and the normal form of the curve. 
We also give characterization of canal surfaces in terms of curves in the set of circles. 
}
\end{abstract}



\section{Introduction}
A space curve is determined by the arc-length, curvature, and torsion up to motions of $\RR^3$. 
In a M\"obius geometric framework, it was shown by Fialkow (\cite{fialkow}) that a space curve is determined by three conformal invariants, conformal arc-length, conformal curvature, and conformal torsion, up to M\"obius transformations. 

The conformal arc-length was found by Liebmann \cite{Li} for space curves and by Pick for plane curves. 
The conformal curvature and torsion were given by Vessiot \cite{Ve}. 
They have been studied by Fialkow using conformal derivations \cite{fialkow}, by Sulanke using Cartan's group theoretical method of moving frames \cite{Su}, and by Cairns, Sharpe, and Webb using normal forms \cite{csw}. 
Conformal torsion was also studied in \cite{MRS} using conformal invariants for pairs of spheres. 

In \cite{La-Oh2} the first and the second authors showed that the conformal arc-length can be considered as ``$\frac12$-dimensional measure'' of the curve of osculating circles. 
This paper is a natural continuation of it. 
We use both curves $\gamma$ of osculating circles and $\sigma$ of osculating spheres, which belong to different spaces, and by doing so we give (hopefully) new formula for the conformal curvature (theorem \ref{characterization_Q_2}). 
We also give moving frames in the Minkowski space using $\gamma$ and $\sigma$. 
Then the conformal curvature and conformal torsion appear in Frenet matrix (this is a just of translation of Sulanke's result \cite{Su} to our context). 
By taking the projection to the Euclidean space which can be realized in $\RR^5_1$ as an affine section of the light cone, we obtain the normal form given in \cite{csw}. 
We give a proof of formulae of conformal curvature and torsion in \cite{csw} using our theorem \ref{characterization_Q_2}. 

Thus, we get geometric and simpler description of integration of preceding studies of conformal invariants of space (and planar) curves. 

We also study canal surfaces and give the characterization of them as curves in the set of circles in $\ss^3$. 

The authors thank Gil Solanes for many useful discussions.

\section{Preliminars}
Let us start with basics in M\"obius geometry which are needed for the study of curves of osculating circles and of osculating spheres. 
\subsection{Realization of $\RR^3$ and $\vector{S}^3$ in Minkowski space $\RR^5_1$}
We start by recalling a commonly used models of a sphere and the Euclidean space in M\"obius  geometry (cf. \cite{Ber, Ce, hertrich}). 
Let us explain in a $3$-dimensional case. 

The {\em Minkowski space} $\RR^5_1$ is $\RR^5$ endowed with an indefinite inner product (the {\em Minkowski product} or the {\em Lorentz form}) given by 
$$ \langle{x},{y}\rangle =-x_0 y_0 +x_1 y_1+x_2 y_2+x_3 y_3+x_4 y_4.$$
The {\em light cone} $\mathcal{C}$ is given by $\mathcal{C}=\{x\in\RR^5_1\,|\,\langle{x},{x}\rangle=0 \}$. 
A vector subspace $W\ne\ RR^5_1$ of $\RR^5_1$  is said to be {\em time-like} if it contains a non-zero time-like vector. 
When $\dim W\ge2$, $W$ is time-like if and only if it intersects the light cone transversally. 

A $3$-sphere $S^3$ or $\RR^3\cup\{\infty\}$ can be identified with the projectivization of the light cone. 
In fact, they can be {\sl isometrically} embedded in $\RR^5_1$ as the intersection of the light cone and a codimension $1$ affine subspace $H$. 
When $H$ can be expressed as $H=\{x\,|\,\esc{x}{n}=-1\}$ for some unit time-like vector $n$, in which case $H$ is tangent to the hyperboloid $\{y\,|\,\esc{y}{y}=-1\}$ at point $n$, $H\cap\mathcal{C}$ is a unit sphere. 
When $H$ can be expressed as $H=\{x\,|\,\esc{x}{n}=-1\}$ for some light-like vector $n$, in which case $H$ is parallel to a codimension $1$ subspace which is tangent to the light cone in $\spa(n)$, $H\cap\mathcal{C}$ is {a paraboloid $\EE^3$. The metric induced from the Lorentz quadratic form on $\EE^3$ is} Euclidean. 
For example, when $n=\big(1/\sqrt2, 1/\sqrt2, 0,0,0\big)$, the Euclidean space can be isometrically embedded as  
\[\RR^3\ni\vector x\mapsto\left(\frac1{\sqrt2}+\frac{\vector x\cdot \vector x}{2\sqrt2},\,-\frac1{\sqrt2}+\frac{\vector x\cdot \vector x}{2\sqrt2},\, \vector x\right)\in\EE^3
=\mathcal{C}\cap\{x\,|\,\esc{x}{n}=-1\}
\,.\]
We use the notation $\ss^3$ and $\EE^3$ to emphasize that they are embedded in $\RR^5_1$. 

\subsection{de Sitter space as the set of codimension $\vector 1$ spheres}
\begin{figure}[htb]
\begin{center}
\psfrag{ssx}{$\Sigma$} \psfrag{lxpur}{ $\sigma^{\bot}$}
\psfrag{xxx}{ $\sigma$} \psfrag{lxx}{ $\spa({\sigma})$}
\psfrag{lamb}{$\Lambda^4$} \psfrag{light}{ light cone}
\psfrag{sinf}{ $\mathbb{S}^3$}
\includegraphics[width=7cm]{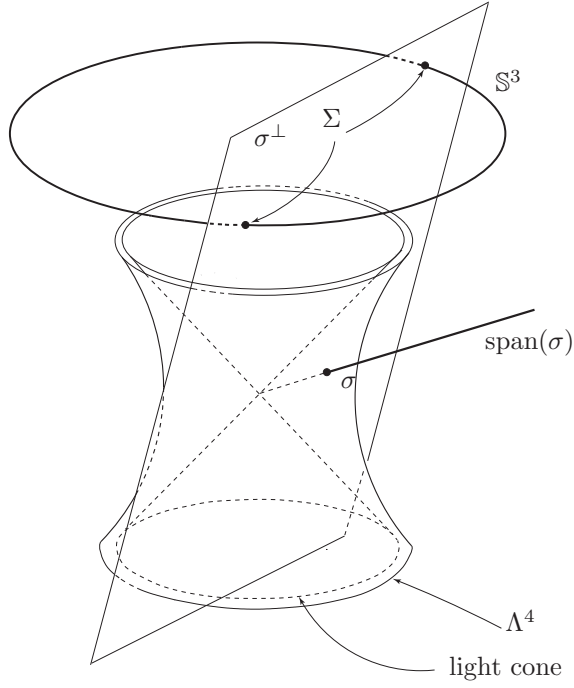}
\end{center}
\caption{The correspondence between de Sitter space and the set of oriented $2$-spheres. \label{fig1}}
\end{figure}

An oriented sphere $\Sigma$ in $\ss^3$ can be obtained as the intersection of an oriented time-like $4$-dimensional subspace of $\RR^5_1$. 
Therefore the set $\mathcal{S}(2,3)$ of oriented spheres in $\ss^3$ can be identified with the Grassmann manifold $\widetilde{G}_{4,5}^{\,-}$ of oriented time-like $4$ dimensional subspaces of $\RR^5_1$. 
By taking the orthogonal complement, we obtain a bijection between $\widetilde{G}_{4,5}^{\,-}$ and the set of oriented space-like lines, which can be identified with the quadric $\Lambda^4=\{x\in\RR^5_1\,|\,\langle{x},{x}\rangle=1\}$ called {\em de Sitter space} (Figure \ref{fig1}). 
The bijection from $\Lambda^4$ to $\mathcal{S}(2,3)$ is given by 
\[\Lambda^4\ni\sigma\mapsto\Sigma=\ss^3\cap(\spa(\sigma))^{\perp}\in\mathcal{S}(2,3)\,,\]
oriented as the boundary of the ball $\ss^3 \cap \{\langle \sigma,.\rangle \leq 0\} \subset \RR^5_1$,
and its inverse is given by 
%
%
\[\mathcal{S}(2,3)\ni\ss^3\cap\spa(u^1,u^2,u^3,u^4)\mapsto\frac{u^1\times u^2\times u^3\times u^4}{\|u^1\times u^2\times u^3\times u^4\|}\in\Lambda^4\, ,\]
where $w=u^1\times u^2\times u^3\times u^4$ is the {\em Lorentz vector product} with respect to $\esc{\cdot}{\cdot}$ that is characterized by $\langle w, u_i\rangle=0$, the norm  $\|w\|=\sqrt{|\esc{w}{w}|}$ being equal to the volume of parallelepiped spanned by $u_i$ which is given by $\sqrt{|\det(\langle u_i, u_j\rangle)|}$, and $\det(w,u^1,u^2,u^3,u^4)>0$ (\cite{La-Oh}\footnote{In fact this gives the opposite sign to \cite{La-Oh}.}). Then $w$ is given by $w_0=-\det(e^0,u^1,u^2,u^3,u^4)$ and $w_i=\det(e^i,u^1,u^2,u^3,u^4)$ $(i\ne0)$. Direct computation shows $(W\ne\RR^5_1)$
\begin{equation}\label{f_<>_vect_prod}
\langle u^1\times u^2\times u^3\times u^4, v^1\times v^2\times v^3\times v^4 \rangle = -\det(\langle u^i, v^j\rangle)
\end{equation}

%
%
%

Similarly, the set of oriented circles in $\ss^2$ can be identified with de Sitter space $\Lambda^3$ in $\RR^4_1$. 

\subsection{Pseudo-Riemannian structure of indefinite Grassmann manifolds} \label{pseudo_riema}
In general, the set of oriented $k$-dimensional subspheres in $\ss^n$ can be identified with the Grassmann manifold $\widetilde{G}_{k+2,n+2}^{\,-}$ of oriented time-like $(k+2)$-dimensional subspaces in $\RR^{n+2}_1$. 
By taking the orthogonal complement we obtain a bijection from $\widetilde{G}_{k+2,n+2}^{\,-}$ to the Grassmann manifold $\widetilde{G}_{n-k,n+2}^{\,+}$ of oriented space-like $(n-k)$-dimensional subspaces in $\RR^{n+2}_1$. 
Thus we need the pseudo-Riemannian structure of these Grassmannians. 

The canonical extension of positive definite inner product to Grassmann algebras is given by 
$$
\langle u^1\wedge \cdots \wedge u^q, v^1\wedge \cdots \wedge v^q\rangle=\det\big(\langle u^i, v^j\rangle\big).
$$
But in our case, as we start with the Minkowski space, it does not fit \eqref{f_<>_vect_prod}. Therefore we agree (after \cite{hertrich} page 280) that the pseudo-Riemannian structures on our Grassmannians are given by  
\begin{equation}\label{f_<>_wedge_prod}
\langle u^1\wedge \cdots \wedge u^q, v^1\wedge \cdots \wedge v^q\rangle=
\left\{
\begin{array}{ll}
\displaystyle -\det\big(\langle u^i, v^j\rangle\big) & \>\>\mbox{if $\mbox{span}(u^i)$ and $\mbox{span}(v^i)$ are time-like},\\[1mm]
\displaystyle \phantom{-}\det\big(\langle u^i, v^j\rangle\big) & \>\>\mbox{if $\mbox{span}(u^i)$ and $\mbox{span}(v^i)$ are space-like}.
\end{array}
\right.
\end{equation}
%

\subsection{Two Grassmannians as the set of oriented circles in $\ss^3$}
An oriented circle $\Gamma$ in $\ss^3$ can be obtained as the intersection of an oriented time-like $3$-dimensional subspace of $\RR^5_1$. 
Therefore the set $\mathcal{S}(1,3)$ of oriented circles in $\ss^3$ can be identified with the Grassmann manifold $\widetilde{G}_{3,5}^{\,-}$ of oriented time-like $3$ dimensional subspaces of $\RR^5_1$. 
By taking the orthogonal complement, we obtain the bijection between $\widetilde{G}_{3,5}^{\,-}$ and the Grassmann manifold $\widetilde{G}_{2,5}^{\,+}$ of oriented space-like $2$ dimensional subspaces of $\RR^5_1$. 


Recall that the wedge product $u\wedge v$ of two vectors in $\RR^5_1$, $u=(u_0, u_1, \ldots, u_4)$ and $v=(v_0, v_1, \ldots, v_4)$, is given by $u\wedge v=(p_{ij})_{0\le i<j\le 4}\in \, \stackrel{2}{\bigwedge}\RR^5_1\cong\RR^{10}$, where the Pl\"ucker coordinates $p_{ij}$ are given by  \setlength\arraycolsep{1pt}
\begin{equation}\label{Plucker_coordinates}
p_{ij}=\left|\begin{array}{cc}
\,u_{i}\, & \,u_{j}\, \\
v_{i} & v_{j} 
\end{array}\right|, 
\end{equation}
and a vectore $\vector{p}=(p_{ij})_{0\le i<j\le 4}$ in $\stackrel{2}{\bigwedge}\RR^5_1\cong\RR^{10}$ is a pure $2$-vector, i.e. the wedge product of two vectors in $\RR^5_1$ if and only if $p_{ij}$ satisfy the Pl\"ucker relation: 
\begin{eqnarray}
p_{01}p_{23}-p_{02}p_{13}+p_{03}p_{12}&=&0,\label{Plucker_relation_1}\\
p_{01}p_{24}-p_{02}p_{14}+p_{04}p_{12}&=&0,\label{Plucker_relation_2}\\
p_{01}p_{34}-p_{03}p_{14}+p_{04}p_{13}&=&0,\label{Plucker_relation_3}\\
p_{02}p_{34}-p_{03}p_{24}+p_{04}p_{23}&=&0,\label{Plucker_relation_4}\\
p_{12}p_{34}-p_{13}p_{24}+p_{14}p_{23}&=&0.\label{Plucker_relation_5}
\end{eqnarray} \setlength\arraycolsep{5pt}
We remark that all of them are not independent. 
For example, the relations (\ref{Plucker_relation_4}) and (\ref{Plucker_relation_5}) can be derived from the rest if $p_{01}\ne0$. 

Then $\widetilde{G}_{2,5}^{\,+}$ can be identified with the set of unit space-like pure $2$-vectors: 
%
\[\widetilde{G}_{2,5}^{\,+} \cong \left\{\vector{p}=(p_{ij})_{0\le i<j\le 4}\,\left|\,\esc{\vector{p}}{\vector{p}}=\sum_{1\le i<j\le 4}p_{ij}{}^2-\sum_{k=1}^{4}p_{0k}{}^2=1, \mbox{ $p_{ij}$ satisfy the Pl\"ucker relation}\right.\right\}.\]
%
%
%
\setlength\arraycolsep{3pt}
%
%
Now the identification between $\mathcal{S}(1,3)$ and $\widetilde{G}_{2,5}^{\,+}$ can be explicitly given by 
\begin{equation}\label{S13<->G25}
\mathcal{S}(1,3)\ni\Gamma=\ss^3\cap(\spa(u,v))^{\perp}\mapsto\frac{u\wedge v}{\|u\wedge v\|}\in \widetilde{G}_{2,5}^{\,+}\,.
\end{equation}
%

\subsection{How to express osculating circles and osculating spheres}
Let $m=m(s)$ be a point in a curve $C$ in $\EE^3\subset\RR^5_1$, where $s$ is the arc-length parameter. 
We always assume that the differential of $m$ never vanishes in what follows. 
Let $\Pi$ be a time-like vector subspace of $\RR^5_1$ of dimension $3$ (or $4$). 
Then the curve $C$ has contact of order $\geq k$ with the circle (or the sphere respectively) $\Pi\cap\EE^3$ if and only if $m(s),m'(s),\cdots,m^{(k)}(s)$ belong to $\Pi$ (remark that the curve is not in $\RR^3$ but in $\EE^3$). 
Therefore, an osculating circle to $C$ at $m(s)$ is given by $\EE^3\cap\spa(m(s),m'(s),m''(s))$, and an osculating sphere is generically given by $\EE^3\cap\spa(m(s),m'(s),m''(s),m'''(s))$. 

In what follows we assume that the curve $C$ is {\em vertex free}, i.e. the osculating circles have contact of order exactly equal to $2$. 
Then the point $\sigma(s)$ in $\Lambda^4$ that corresponds to the osculating sphere at point $m(s)$ is given by 
\begin{equation}\label{f_sigma}
\sigma(s)=\frac{m(s)\times m'(s)\times m''(s)\times m'''(s)}{\|m(s)\times m'(s)\times m''(s)\times m'''(s)\|}\,.
\end{equation}
Note that $\sigma'(s)$ can be expressed as 
\begin{equation}\label{f_sigma'}
\sigma'(s)=m(s)\times m'(s)\times m''(s)\times \big(a\,m'''(s)+b\,m^{(4)}(s)\big)
\end{equation}
for some $a,b\in\RR$. 
Let us further assume that $\|\sigma'(s)\|$ never vanishes. 
As $\esc{\sigma}{\sigma'}=0$ we have $\dim\spa(\sigma,\sigma')=2$. 
The formulae \eqref{f_sigma} and \eqref{f_sigma'} imply 
$$\spa(\sigma(s),\sigma'(s))=\big(\spa(m(s),m'(s),m''(s))\big)^{\perp}.$$
Therefore, the osculating circle $\Gamma(s)$ to $C$ at point $m(s)$ corresponds to a point $\gamma(s)$ in $\widetilde{G}_{2,5}^{\,+}$ given by 
\[
\gamma(s)=\pm\frac{\sigma(s)\wedge\sigma'(s)}{\|\sigma(s)\wedge\sigma'(s)\|}\,.
\]
For simplicity's sake, we assume that the sign is $+$ in what follows. 
If we denote the derivative by the arc-length parameter $l$ of the curve of osculating spheres $\sigma$ in $\Lambda^4$  by putting $\Pt{\phantom{\sigma}}$ above, we have 
\begin{equation}\label{gamma=sigma_wedge_sigma_dot}
\gamma=\sigma\wedge\Pt\sigma
\end{equation}

\medskip
An osculating circle to a curve $C$ in $\EE^2$ (in this case, the osculating sphere is constantly equal to $\EE^2$) can be expressed in a similar way as \eqref{f_sigma}. 
In this case, as $\|m(s)\times m'(s)\times m''(s)\|=1$ (\cite{La-Oh2}), a point $\gamma$ in $\Lambda^3$ which corresponds to the osculating circle is given by $\gamma=m(s)\times m'(s)\times m''(s)$. 

\smallskip
Similarly, an osculating circle to a curve in $\EE^3$ can be expressed by %
\begin{equation}\label{f_gamma_XX}
\gamma(s)=m(s)\wedge m'(s)\wedge m''(s) \in \widetilde{G}_{3,5}^{\,-}. 
\end{equation}

\subsection{Our notations and assumptions}
We express a curve in $\EE^2$ or $\EE^3$ by $C$ and a point on it by $m$. 
We always assume that $C$ is vertex-free. 
The letters $\gamma$ and $\sigma$ express the osculating circle and sphere respectively unless otherwise mentioned. 
There are one or two natural conformally invariant parameters: one is the conformal arc-length which will be denoted by $\rho$, and the other is the arc-length parameter of the curve of osculating spheres which will be denoted by $l$. 
The latter appears only when $C$ is a space {curve}. 
We also use the arc-length $s$ of $C$. 
In order to avoid confusion, we use different notations to express the derivatives by these parameters. 
The derivatives by $s, \rho$, and $l$ are expressed by putting ${}', \c{\phantom{\sigma}}$, and $\Pt{\phantom{\sigma}}$ respectively. 
For simplicity's sake, we put a technical assumption that $\frac{dl}{d\rho}$ never vanishes when $C$ is a space curve. 

\subsection{Null curve of osculating circles and conformal arc-length}
Let us first consider a smooth one parameter family of circles in $\RR^2$ given by their centers $\omega(t)$ and radii $r(t)$. 
Locally speaking, these circles may admit an envelope formed of two curves, or exceptionally one curve; this is the case we will consider here, because we consider the family of osculating circles to a curve.   
\begin{figure}[ht]\label{osc_circ}
\begin{center}
\includegraphics[width=.4\linewidth]{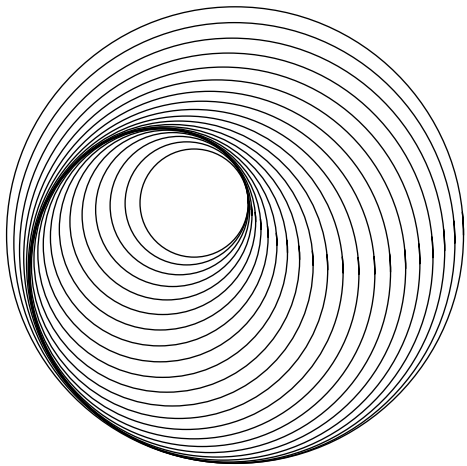}%
\caption{Osculating circles of a plane curve with monotone
curvature\label{oscul}}
\end{center}
\end{figure}

Then $\|\omega'(t)\|=|r'(t)|$ everywhere. 
The fact that the osculating circles of a planar arc with monotone curvature are nested along the arc as is illustrated in Figure \ref{osc_circ} was observed in the beginning of $20$th century by Kneser (\cite{kneser}).

In our language, 
\begin{lemm}\label{lemma_nullity_osc_circles}{\rm (\cite{La-Oh2})}
A curve $\gamma\subset\Lambda^3$ of osculating circles to a curve $C$ in $\EE^2$ is light-like, and the point $m(s)$ on $C$ can be given by $m(s)=\EE^2\cap\spa(\gamma'(s))$. 

Similarly, when $C$ is a curve in $\EE^3$, a curve $\gamma$ of osculating circles is a null curve in $\widetilde{G}_{2,5}^{\,+}$ with $\gamma'$ being a pure vector. 
If $\Pi(s)$ is a plane in $\RR^5_1$ which corresponds to $\gamma'(s)\in\stackrel{2}{\bigwedge}\RR^5_1$, then it is tangent to the light cone in a line $\spa(m(s))$. 
\end{lemm}

%
%
\begin{prop}\label{prop_gamma_circcirc}{\rm (\cite{La-Oh2})} 
Let $\gamma$ denote a curve of osculating circles to a curve $C$ in $\EE^2$ or $\EE^3$. 
Then the $1$-form $\sqrt[4]{\big\langle\frac{d^2\gamma}{dt^2},\frac{d^2\gamma}{dt^2}\big\rangle}\,dt$ is independent of the parameter $t$. 
Let $\rho$ be a parameter so that the $d\rho=\sqrt[4]{\big\langle\frac{d^2\gamma}{dt^2},\frac{d^2\gamma}{dt^2}\big\rangle}\,dt$. 
In other words, the parameter $\rho$ can be characterized by 
\begin{equation}\label{gamma_circ_circ}
\esc{\cc\gamma}{\cc\gamma}=1,
\end{equation}
where $\cc\gamma=d^2\gamma/d\rho^2$. 
It can be uniquely determined up to $\rho\mapsto\pm\rho+c$ for some constant $c$. 
\end{prop}

This parameter is called the {\em conformal arc-length} of $C$. 
Our assumption that the curve $C$ is vertex-free guarantees that $\rho$ serves as a non-singular parameter of $C$. 

\section{Curves in $\EE^2$ or $\ss^2$}

Let $\gamma\subset\Lambda^3$ be an osculating circle to a curve $C$ in $\EE^2$ at a point $m$. 
Then they form a curve in de Sitter space $\Lambda^3$. 
We express the derivative by the conformal arc-length $\rho$ by putting $\c{\phantom{\sigma}}$ above. 

\subsection{The moving frame and Frenet formula}
Our moving frames consist of two space-like vectors and two light-like vectors, instead of three space-like vectors and a time-like vector, since we take a light-like vector in $\spa(m)$ as our first frame, which we denote by $n$. 
Another light-like vector, which comes last in our frames, is chosen so that $\esc{n}{n^{\ast}}=-1$. The middle two space-like vectors of our frames are taken from  an orthonormal basis of $(\spa(n,n^{\ast}))^{\perp}$. 
The moving frames of this form are called {\em isotropic orthonormal} frames. 

Let us choose the first vector $n$ of our moving frames in $\spa(m)$ so that $\|\c n\|=1$. 
Then we can take $n=\c\gamma$. 
Our second vector $v_1$ is $\c n$. 
Then the $x$-coordinate of the normal form can be given by $\esc{m(\rho)}{v_1}$. 
We choose our third vector $v_2$ in $(\spa(m))^{\perp}$ so that $\esc{v_i}{v_j}=\delta_{ij}$ and that  the $y$ coordinate of the normal form can be given by $\esc{m(\rho)}{v_2}$. 
The $x$-axis of the normal form should be the osculating circle. 
Therefore, we can take $v_2=\gamma$. 
Our last vector $n^{\ast}$ is a light-like vector in $(\spa(v_1,v_2))^{\perp}$ that satisfies $\esc{n}{n^{\ast}}=-1$. 

Put (cf. \cite{Su})
\begin{equation}\label{f_Q_plane}
Q_2=-\frac12\big\langle\cc{n}, \cc{n}\big\rangle\,. 
\end{equation}

Then $\big\langle\ccc\gamma, \ccc\gamma\big\rangle=-2Q_2$. 
Lemma \ref{lemma_nullity_osc_circles} and proposition \ref{prop_gamma_circcirc} imply that the $\Big\langle \mbox{\large $\frac{d^{\,i}{\gamma}}{d\rho^i}, \frac{d^{\,j}{\gamma}}{d\rho^j} $}\Big\rangle$ for small $i,j$ are given by table \ref{table_gammaij_circ}. 
It implies that our last frame is given by $n^{\ast}=-Q_2\c\gamma+\ccc\gamma=-Q_2\,n+\cc{n}$. 
\begin{table}[htb]
\begin{center}
\par\noindent
\begin{tabular}{|c|c|c|c|c|c|c|c|}
\hline
$\langle\cdot,\cdot\rangle$ $\phantom{\stackrel{{a}}{0}}$
& $\gamma$  &  $\c\gamma$  &  $\cc\gamma$  & $\ccc\gamma$ & $\gamma^{(4)}$ & $\gamma^{(5)}$ & $\gamma^{(6)}$\\[1mm] 
\hline
 $\!\!\phantom{\stackrel{{a}}{0}}\!\!\gamma$ & $1$ & $0$ & $0$ & $0$ & $1$ & $0$ & $2Q_2$ \\[1mm] 
\hline
 $\!\!\phantom{\stackrel{{a}}{0}}\!\!\c\gamma$ & $0$ & $0$ & $0$ & $-1$ & $0$ & $-2Q_2$ &  \\[1mm] 
\hline
 $\!\!\phantom{\stackrel{{a}}{0}}\!\!\cc\gamma$ & $0$ & $0$ & $1$ & $0$ & $2Q_2$ &  &  \\[1mm] 
\hline
 $\!\!\phantom{\stackrel{{a}}{0}}\!\!\ccc\gamma$ & $0$ & $-1$ & $0$ & $-2Q_2$ &  &  &  \\[1mm]%
\hline 
\end{tabular}
\end{center}
\caption{A table of $\Big\langle \mbox{\large $\frac{d^{\,i}{\gamma}}{d\rho^i}, \frac{d^{\,j}{\gamma}}{d\rho^j} $}\Big\rangle$} 
\label{table_gammaij_circ}
\end{table}

We remark that our moving frames 
\[n=\c\gamma, \ v_1=\c{n}=\cc\gamma, \ v_2=\gamma, \ \ n^{\ast}=-Q_2\, n+\cc{n}=-Q_2\c\gamma+\ccc\gamma\] 
also serve as moving frames of the curve $\gamma\subset\Lambda^3$ of osculating circles 
(Remark that $n=\c\gamma$ is not necessarily equal to a point $m$ in the curve $C\subset\EE^2$). 

\medskip
The Frenet formula with respect to our moving frames is given by 
\begin{equation}\label{Frenet_matrix_plane}
\frac{d}{d\rho}
\left(\begin{array}{l}
n \\ v_1 \\ v_2 \\ n^{\ast}
\end{array}\right)
=
\left( \begin{array}{cccc}
0&  1&  0&  0 \\
Q_2&  0& 0 & 1 \\
1&  0&  0& 0 \\
0&  Q_2 &  1&  0
\end{array}\right)
\left(\begin{array}{l}
n \\ v_1 \\ v_2 \\ n^{\ast}
\end{array}\right).
\end{equation}

See remark \ref{remark_other_refe}. 

\subsection{Normal form} 
The normal form of a planar curve can be obtained from that of a space curve (see subsection \ref{subsec_normal_form_space}) by putting $T=0$ and $Q=Q_2$ and forgetting the $z$-coordinate. 
If we only look for the normal form of a plane curve, the computation is simpler. 
Thanks to table \ref{table_gammaij_circ}, we do not have to rewrite $n^{(i)}=\gamma^{(i+1)}$ in terms of the frames. 
The direction vectors of $x$- and $y$- axes at $\rho=0$ are given by $\cc\gamma(0)$ and $\gamma(0)$ respectively. 

\section{Curves in $\EE^3$ or $\ss^3$}
%
\subsection{Osculating spheres and conformal torsion}
Let $m$ be a point on a curve $C$ in $\EE^3$. 
Let $\gamma\subset\widetilde{G}_{2,5}^{\,+}$ be a curve of osculating circles, and $\sigma\subset\Lambda^4$ a curve of osculating spheres. 
Since an osculating sphere intersects an infinitesimally close osculating sphere in an osculating circle, $\sigma$ is a space-like curve in $\Lambda^4$. 
Let $\rho$ be the conformal arc-length and $l$ the arc-length parameter of the curve $\sigma$. 
We express the derivatives by $l$ and $\rho$ by putting above $\Pt{\phantom{\sigma}}$ and $\c{\phantom{\sigma}}$ respectively. 
Put 
\begin{equation}\label{f_T}
T=\frac{dl}{d\rho}\,.
\end{equation}
As $dl$ measures the infinitesimal angle variation, $T$ measures how an osculating sphere rotates around an osculating circle with respect to the conformal arc-length. 
For simplicity's sake, let us assume $T>0$ in what follows. 
It is proved in \cite{Ro-Sa} that $T$ coincides with the {\em conformal torsion} up to sign. 
As is pointed out in \cite{Sh}, the conformal torsion can be determined up to sign. 
We remark that the conformal torsion is identically equal to $0$ if and only if $C$ is a planar or a spherical curve. 
\begin{lemm}\label{coro_sigma_dotdot} 
We have $\esc{\Ppt\sigma}{\Ppt\sigma}=1$. 
\end{lemm}

\begin{demo}
Since $\esc{\sigma}{\sigma}=1$ and $\esc{\Pt\sigma}{\Pt\sigma}=1$, we have $\esc{\sigma}{\Pt\sigma}=0$ and $\esc{\sigma}{\Ppt\sigma}=-1$. 
On the other hand, since $\gamma=\sigma\wedge\Pt\sigma$ we have $\Pt\gamma=\sigma\wedge\Ppt\sigma$. 
Therefore, the formula \eqref{f_<>_wedge_prod} implies 
\[\esc{\Pt\gamma}{\Pt\gamma}=\left|\begin{array}{cc}
\esc{\sigma}{\sigma} & \esc{\sigma}{\Ppt\sigma} \\
\esc{\sigma}{\Ppt\sigma} & \esc{\Ppt\sigma}{\Ppt\sigma}
\end{array}\right|.\]
As $\esc{\Pt\gamma}{\Pt\gamma}=0$ by lemma \ref{lemma_nullity_osc_circles}, it implies $\esc{\Ppt\sigma}{\Ppt\sigma}=1$. 
\end{demo}

\begin{prop}\label{prop_conf_torsion}
The conformal torsion $T$ satisfies
\[
T=\frac1{\sqrt[4]{\langle\Ppt\gamma,\Ppt\gamma\rangle}}=\frac1{\sqrt[4]{\langle\Pppt\sigma,\Pppt\sigma\rangle-1}}\,.
\]
\end{prop}

\begin{demo}
The first equality comes from proposition \ref{prop_gamma_circcirc}, the fact that $d\rho=\sqrt[4]{\big\langle\frac{d^2\gamma}{dt^2},\frac{d^2\gamma}{dt^2}\big\rangle}\,dt$ for any parameter $t$ (\cite{La-Oh2}). 

Since $\esc{\sigma}{\sigma}=\esc{\Pt\sigma}{\Pt\sigma}=\esc{\Ppt\sigma}{\Ppt\sigma}=1$, we have 
$\esc{\sigma}{\Pt\sigma}=\esc{\Pt\sigma}{\Ppt\sigma}=\esc{\Ppt\sigma}{\Pppt\sigma}=0$ and hence 
$\esc{\sigma}{\Ppt\sigma}=\esc{\Pt\sigma}{\Pppt\sigma}=-1$ and $\esc{\sigma}{\Pppt\sigma}=0$. 
Therefore the formula \eqref{f_<>_wedge_prod} implies 
\[
\esc{\Pt\si\wedge\Ppt\si}{\Pt\si\wedge\Ppt\si}=1, \>
\esc{\Pt\si\wedge\Ppt\si}{\si\wedge\Pppt\si}=-1, \>\mbox{ and }\>
\esc{\si\wedge\Pppt\si}{\si\wedge\Pppt\si}=\esc{\Pppt\si}{\Pppt\si}\,.
\]
As $\gamma=\si\wedge\Pt\si$ and therefore $\Ppt\gamma=\si\wedge\Pppt\si+\Pt\si\wedge\Ppt\si$, 
we have $\esc{\Ppt\gamma}{\Ppt\gamma}=\langle\Pppt\sigma,\Pppt\sigma\rangle-1$, which completes the proof. 
%
%
%
%
%
\end{demo}
It follows that $\langle d^i\sigma/dl^i, d^j\sigma/dl^j\rangle$ for small $i,j$ are given by table \ref{table_sigma_ij}. 
\begin{table}[h]
\begin{center}
\par\noindent
\begin{tabular}{|c|c|c|c|c|c|}
\hline
$\langle\cdot,\cdot\rangle$ $\phantom{\stackrel{{a}}{0}}$
& $\sigma$  &  $\Pt\sigma$  &  $\Ppt\sigma$  & $\Pppt\sigma$ & $\sigma^{(4)}$ \\[1mm] 
\hline
 $\!\!\phantom{\stackrel{{a}}{0}}\!\!\sigma$ & $1$ & $0$ & $-1$ & $0$ & $1$ \\[1mm] 
\hline
 $\!\!\phantom{\stackrel{{a}}{0}}\!\!\Pt\sigma$ & $0$ & $1$ & $0$ & $-1$ & $0$ \\[1mm] 
\hline
 $\!\!\phantom{\stackrel{{a}}{0}}\!\!\Ppt\sigma$ & $-1$ & $0$ & $1$ & $0$ & $-(1+T^{-4})$ \\[1mm] 
\hline
 $\!\!\phantom{\stackrel{{b}}{0}}\!\!\Pppt\sigma$ & $0$ & $-1$ & $0$ & $1+T^{-4}$ & $-2T^{-5}\Pt T$ \\[1mm] 
\hline 
\end{tabular}
\end{center}
\caption{A table of $\Big\langle \mbox{\large $\frac{d^{\,i}\sigma}{{dl}^{i}}, \frac{d^{j}\sigma}{{dl}^{j}} $}\Big\rangle$}
%
\label{table_sigma_ij}
\end{table}

\begin{lemm}\label{lemm_lin_indep}
Under our assumption, i.e. if the curve $C$ is vertex-free and $\frac{dl}{d\rho}$ never vanishes, the five vectors $\sigma, \Pt\sigma, \Ppt\sigma, \Pppt\sigma$, and $\sigma^{(4)}$ are linearly independent. 
\end{lemm}

\begin{demo}
Suppose $a_0\si+a_1\Pt\si+a_2\Ppt\si+a_3\Pppt\si+a_4\si^{(4)}=\vector{0}$ for some $a_0,\cdots a_4\in\RR$. 
Put $G_3=\langle\Pppt\si,\Pppt\si\rangle$ and $G_4=\langle\si^{(4)}, \si^{(4)}\rangle$. 
Then, by taking pseudo inner product with $\sigma, \cdots, \si^{(4)}$, we have 
\[
\left(\begin{array}{ccccc}
1 & 0 & -1 & 0 & 1 \\
0 & 1 & 0 & -1 & 0 \\
-1 & 0 & 1 & 0 & -G_3 \\
0 & -1 & 0 & G_3 & \frac12\Pt{G_3} \\
1 & 0 & -G_3 & \frac12\Pt{G_3} & G_4
\end{array}
\right)
\left(\begin{array}{c}
a_0\\a_1\\a_2\\a_3\\a_4
\end{array}
\right)
=\vector{0}\,.
\]
The determinant of the coefficient matrix is equal to $(1-G_3)^3=-T^{-12}\ne0$. 
\end{demo}

\subsection{Moving frame in $\RR^5_1$ and Frenet formula}
\begin{prop}{\rm (\cite{yun})} 
A point $m$ on the curve $C$ can be expressed in terms of the osculating spheres as 
\begin{equation}\label{m_sigma+sigma_dotdot}
m=c(\sigma+\Ppt\sigma)\hspace{0.5cm}(c\in\RR\setminus\{0\}). 
\end{equation}
%
\end{prop}

\begin{demo}
Since $\esc{\sigma}{\Ppt\sigma}=-1$ and $\esc{\Ppt\sigma}{\Ppt\sigma}=1$ by lemma \ref{coro_sigma_dotdot} and its proof, $\sigma+\Ppt\sigma$ is a light-like vector. 

On the other hand, as $\sigma$ can be expressed as $\sigma=\varphi \,m\times\Pt m\times \Ppt m\times \Pppt m$ for some function $\varphi$, $\Ppt\sigma$ can be expressed as a linear combination of vectors of the form $m\times \Pt m \times m^{(i)}\times m^{(j)}$. 
Therefore, $\esc{m}{\sigma+\Ppt\sigma}=0$, which completes the proof, as $(\spa(m))^{\perp}$ is tangent to the light cone in $\spa(m)$. 
\end{demo}

We choose the conformal arc-length, not $l$, for the parameter of the Frenet formula, as is the case in most preceding studies. 
Let us choose the first vector $n$ of our moving frames in $\spa(m)$ so that $\|\c n\|=1$. 
Then we can take $n=T(\sigma+\Ppt\sigma)$. 
Our second vector $v_1$ is $\c n$. 
Then the $x$-coordinate of the normal form can be given by $\esc{m(\rho)}{v_1}$. 
We choose our third and fourth vectors $v_2$, and $v_3$ in $(\spa(m))^{\perp}$ so that $\esc{v_i}{v_j}=\delta_{ij}$ and that  the $y$, and $z$ coordinates of the normal form can be given by $\esc{m(\rho)}{v_2}$, and $\esc{m(\rho)}{v_3}$ respectively. 
The $xy$-plane of the normal form should be the osculating sphere. 
Therefore we can take $v_3=-\sigma$. The choice of the frame is important to get the same normal form as in \cite{csw}. We put $-$ here so that the Frenet matrix and the normal form fit with those in \cite{Su} and \cite{csw} respectively. 
A sphere which corresponds to $v_2$ should intersect the osculating sphere orthogonally in the osculating circle. 
Therefore we can take $v_2=\Pt\sigma$. 
Our last vector $n^{\ast}$ is a light-like vector in $(\spa(v_1,v_2,v_3))^{\perp}$ that satisfies $\esc{n}{n^{\ast}}=-1$. 

We remark that our moving frames 
\[n=T(\sigma+\Ppt\sigma), \> v_1=\c n=T\Pt T(\sigma+\Ppt\sigma)+T^2(\Pt\sigma+\Pppt\sigma), \> 
v_2=\Pt\sigma, \> v_3=-\sigma, \> n^{\ast}\]
also serve as moving frames of the curve $\sigma\subset\Lambda^4$ of osculating spheres. 
Note that $n$ is not necessarily equal to a point $m$ in the curve $C\subset\EE^3$. 
In order to have $\|\c n\|=1$, we enlarge or shrink the point vector of a curve in the light cone keeping the parameter to be the conformal arc-length. 

\begin{prop}\label{lem_Frenet_space}{\rm (\cite{Su})}
Put
\begin{equation}\label{def_Q_3}
Q=\big\langle \c{n}, \c{n}{}^{\ast}\big\rangle = -\big\langle{\cc{n}}, {n^{\ast}}\big\rangle\,.
\end{equation}
Then the Frenet formula with respect to our moving frames is given by \setlength\arraycolsep{4pt}
\begin{equation}\label{Frenet_matrix_space}
\frac{d}{d\rho}
\left(\begin{array}{l}
n \\ v_1 \\ v_2 \\ v_3 \\ n^{\ast}
\end{array}\right)
=
\left( \begin{array}{ccccc}
 0 & 1 & 0 & 0 & 0 \\
 Q & 0 & 0 & 0 & 1 \\
 1 & 0 & 0 & T & 0 \\
 0 & 0 & -T & 0 & 0 \\
 0 & Q & 1 & 0 & 0 
\end{array}\right)
\left(\begin{array}{l}
n \\ v_1 \\ v_2 \\ v_3 \\ n^{\ast}
\end{array}\right).
\end{equation}
\end{prop}

\begin{demo}
The formula \eqref{f_T} implies $\c{u}=T\Pt{u}$ for any $u$. 
The first, third, and fourth rows of the matrix follow from the definition of the frames. 
The other two rows can be obtained by derivating the scalar products of the frames, i.e. $\esc{n}{v_1}$ etc. 
Remark that ``$n$- (or $n^{\ast}$-)coordinate'' of a vector $u$ can be given by $-\esc{u}{n^{\ast}}$ (or $-\esc{u}{n}$ respectively). 
\end{demo}

Since $\cc{n}=\c{v_1}=Qn+n^{\ast}$ we have $\big\langle\cc{n}, \cc{n}\big\rangle=\big\langle Qn+n^{\ast}, Qn+n^{\ast}\big\rangle=-2Q$. 
Therefore, 
\begin{equation}\label{f_KQ_space} 
Q=-\frac12\big\langle\cc{n}, \cc{n}\big\rangle\,.
\end{equation}
This quantity $Q$ is called the {\em conformal curvature} (denoted by $\lambda_1$ in \cite{Su}, whereas the conformal torsion $T$ is denoted by $\pm\lambda_2$ in \cite{Su}). 

As $\cc{n}=Qn+n^{\ast}$, our last frame is given by $n^{\ast}=-Q\,n+\cc{n}$ as in the case of planar curve. 
Thus our isotropic orthonormal moving frames are 
\[n\in\spa(m), \, v_1=\c{n}, \, v_2=\frac{\c\sigma}{\|\c\sigma\|}, \, v_3=-\sigma, \, n^{\ast}=-Qn+\cc{n}.\]

\begin{theo}\label{characterization_Q_2}
Let $\gamma\subset\widetilde{G}_{2,5}^{\,+}$ be a curve of osculating circles and $\sigma\subset\Lambda^4$ a curve of osculating spheres. 
Then the conformal curvature $Q$ and the conformal torsion $T$ satisfy 
\begin{eqnarray}
T&=&\displaystyle \sqrt{\big\langle\c\sigma, \c\sigma\big\rangle}\,,\label{f_T2} \\
Q&=&\displaystyle -\frac12\big\langle\ccc\gamma, \ccc\gamma\big\rangle+3\big\langle\c\sigma, \c\sigma\big\rangle. \label{f_Q_gamma_sigma}
\end{eqnarray}
%
\end{theo}

\begin{demo}
The first equation is trivial from the definition \eqref{f_T} of $T$ as $\big\langle\c\sigma, \c\sigma\big\rangle=\big\langle T\Pt\sigma, T\Pt\sigma\big\rangle=T^2$. 

Since $\gamma=\sigma\wedge\Pt\sigma=v_3\wedge v_2$, we have 
\[\begin{array}{rcl}
\c\gamma&=&n\wedge v_3\\
\ccc\gamma&=&-2Tv_1\wedge v_2-\c{T}n\wedge v_2+(Q-T^2)n\wedge v_3-v_3\wedge n^{\ast}\,.
\end{array}\]
It follows that $\big\langle\ccc\gamma, \ccc\gamma\big\rangle=4T^2-2(Q-T^2)=-2Q+6T^2$, which implies the second equation. 
\end{demo}

In the plane curve case, the last frame $n^{\ast}$ was given by $n^{\ast}=\frac12\big\langle\ccc\gamma, \ccc\gamma\big\rangle\,\c\gamma+\ccc\gamma$. 
But in the space curve case, the proof of the above lemma shows that $\frac12\big\langle\ccc\gamma, \ccc\gamma\big\rangle\,\c\gamma+\ccc\gamma$, or in general any vector of the form $c\c\gamma+\ccc\gamma$, is not a pure vector (see section\ref{pseudo_riema}), in other words, although it is a null vector, it does not define a point in $\EE^3$. 

\begin{remark}\rm 
Let us consider the asymptotic behavior of the frames as the curve $C$ degenerates to a plane curve, which corresponds to $T\to0$. 
We have
\[v_2\wedge v_3=\gamma, \> v_1\wedge v_3=\cc\gamma+T^2(\gamma-\Pt\sigma\wedge\Ppt\sigma), \> v_1\wedge v_3=T(\Pt{T}\gamma-\Pt{T}\sigma\wedge\Ppt\sigma-T\sigma\wedge\Pppt\sigma)\]
Therefore, as $T$ goes to $0$, $v_2\wedge v_3$ and $v_1\wedge v_3$ give the $x$- and $y$-axes of $\EE^2$ respectively as is desired. 
\end{remark}

\begin{remark}\label{remark_other_refe}\rm 
Let us give a comment on the correspondence of our expression to some others in references of M\"obius geometry. 
Let us consider  the 5-dimensional Lorentz space equipped with the indefinite form given by 
$\langle{\xi},{\eta}\rangle=\xi_1\eta_1+\xi_2\eta_2+\xi_3\eta_3-\xi_0\eta_4-\xi_4\eta_0.$ 
Then the Euclidean space can be realized in this space as $\EE^3_L=\{(0,x,y,z,0)\,|\,(x,y,z)\in\RR^2\}$. 
Let $\mathcal{M}\mbox{\it \"ob}_3$ be the M\"obius group. 
The moving frames along a curve $C$ define a map $g=(n\,v_1\,v_2\,v_3\,n^{\ast})$ from $C$ to $\mathcal{M}\mbox{\it \"ob}_3$, where the vectors are considered as column vectors. 
Then, our matrix in the Frenet formula is a transposed matrix of $g^{-1}\frac{dg}{d\rho}$. 
We remark that there is a difference of a choice of signature in \cite{csw}, where the indefinite form is given by $\langle{\xi},{\eta}\rangle=\xi_1\eta_1+\xi_2\eta_2+\xi_3\eta_3+\xi_0\eta_4+\xi_4\eta_0$, causing a difference in sign of the last row and column of the matrix. 
\end{remark}

\subsection{Normal form}\label{subsec_normal_form_space}
The point $n=n(\rho)$ does not necessarily stay on some constant Euclidean model $\EE^3$. 
In order to have the normal form as in \cite{csw} we have to first project a point $n(\rho)$ to a point $m(\rho)$ in $\EE^3$ by a radial projection from the origin, and then use a M\"obius transformation of $\EE^3\cup\{\infty\}$ which sends $m^{\ast}(\rho_0)=\EE^3\cap\spa(n^{\ast}(\rho_0))$ to $\infty$. 
We may assume without loss of generality that $\rho_0=0$ and that our Euclidean model $\EE^3$ is given by 
$\EE^3=\{v\in\RR^5_1\,|\,\esc{v}{n^{\ast}(0)}=-1\}.$ 
Let $f(\rho)$ be a function such that $m(\rho)$ is given by $m(\rho)=f(\rho)n(\rho)$. 
Remark that $m^{\ast}$ is on the osculating circle to $C$ as $\esc{\sigma}{n^{\ast}}=\esc{\Pt\sigma}{n^{\ast}}=0$. 
We use the Frenet formula \eqref{Frenet_matrix_space} to express $n^{(i)}(\rho)$ in terms of the frames. 
Then 
\[\begin{array}{rcl}
-1=\esc{f(\rho)n(\rho)}{n^{\ast}(0)}
&=&\displaystyle f(\rho)\Big\langle n(0)+\c n(0)\,\rho+\frac{\cc n(0)}{2!}\,\rho^2+O(\rho^3), n^{\ast}(0)\Big\rangle \\[2mm]
&=&\displaystyle f(\rho)\Big(-1-\frac{Q(0)}2\,\rho^2+O(\rho^3)\Big),
\end{array}\]
which implies that 
\[f(\rho)=1-\frac{Q(0)}2\,\rho^2+O(\rho^3),\]
and therefore, the normal form is given by 
\[\begin{array}{rcl}
x&=&\displaystyle \big\langle{m(\rho)}, {v_1(0)}\big\rangle
=f(\rho)\,\big\langle{n(\rho)}, {v_1(0)}\big\rangle\\[2mm]
&=&\displaystyle \Big(1-\frac{Q(0)}2\,\rho^2+O(\rho^3)\Big)\Big(\rho+\frac{Q(0)}3\,\rho^3+O(\rho^4)\Big)\\[3mm]
&=&\displaystyle \rho-\frac{Q(0)}6\,\rho^3+O(\rho^4),\\[3mm]
y&=&\displaystyle \big\langle{m(\rho)}, {v_2(0)}\big\rangle
=f(\rho)\,\big\langle{n(\rho)}, {v_2(0)}\big\rangle\\[2mm]
&=&\displaystyle \Big(1-\frac{Q(0)}2\,\rho^2+O(\rho^3)\Big)\Big(\frac1{3!}\,\rho^3+\frac{2Q(0)-T(0)^2}{5!}\,\rho^5+O(\rho^6)\Big)\\[4mm]
&=&\displaystyle \frac{1}{3!}\,\rho^3+\left(\frac{2Q(0)-T(0)^2}{5!}-\frac{Q(0)}{2\cdot3!}\right)\,\rho^5+O(\rho^6)\,, \\[4mm]
z&=&\displaystyle \big\langle{m(\rho)}, {v_3(0)}\big\rangle
=f(\rho)\,\big\langle{n(\rho)}, {v_3(0)}\big\rangle\\[1mm]
&=&\displaystyle \Big(1-\frac{Q(0)}2\,\rho^2+O(\rho^3)\Big)\Big(\frac{T}{4!}\,\rho^4+\frac{\c T}{5!}\,\rho^5+O(\rho^6)\Big)\\[3mm]
&=&\displaystyle \frac{T}{4!}\,\rho^4+\frac{\c T}{5!}\,\rho^5+O(\rho^6)\,.
\end{array}\]
Since 
\[\frac{1}{3!}\,x^3=\frac{1}{3!}\,\rho^3-\frac{Q(0)}{2\cdot3!}\,\rho^5+O(\rho^6)\]
the normal form is given by 
\begin{equation}\label{f_normal_form_space}
\begin{array}{rcl}
y&=&\displaystyle \frac{1}{3!}\,x^3+\frac{(2Q-T^2)}{5!}\,x^5+O(x^6),\\[2mm]
z&=&\displaystyle \frac{T}{4!}\,x^4+\frac{\c{T}}{5!}\,x^5+O(x^6),
\end{array}\end{equation}
\noindent as in \cite{csw}.
\subsection{Moving frames for $\vector{\gamma}$ in $\stackrel{2}{\bigwedge}\RR^5_1\cong\RR^{10}_4$}
Our moving frames $n, v_1, v_2, v_3, n^{\ast}$ produce moving frames of $\stackrel{2}{\bigwedge}\RR^5_1\cong\RR^{10}_4$ in turn, which consist of $4$ space-like vectors 
\[n\wedge n^{\ast},\, v_1\wedge v_2,\, v_1\wedge v_3,\, v_2\wedge v_3=\gamma, \]
and $6$ null vectors
\[n\wedge v_1,\, n\wedge v_2,\, n\wedge v_3=\c\gamma,\, v_1\wedge n^{\ast},\, v_2\wedge n^{\ast},\, v_3\wedge n^{\ast}.\]
The Frenet formula for them can be obtained using \eqref{f_<>_wedge_prod} and \eqref{Frenet_matrix_space}.

\section{Euclidean expression of conformal invariants}
Let us express the conformal invariants of curves in terms of derivatives of curvature $\kappa$ and torsion $\tau$ of a curve $C\subset\EE^3$ with respect to the arc-length $s$. 

Put $\nu=\sqrt{(\kappa')^2+\kappa^2\tau^2}$, where ${}'$ means $\frac{d}{ds}$. 
It never vanishes if and only inf $C$ is vertex-free. 
Then the conformal arc-length $\rho$ satisfies (\cite{Ta})
\begin{equation}\label{f_drho_ds}
d\rho=\sqrt{\nu}\,ds=\sqrt[4]{(\kappa')^2+\kappa^2\tau^2}\,ds.
\end{equation}
\begin{prop}{\rm (\cite{csw})}
The conformal curvature $Q$ and conformal torsion $T$ can be expressed as 
\begin{eqnarray}
Q&=&\displaystyle \frac{4(\nu\,''-\kappa^2\nu)\nu-5(\nu\,')^2}{8\nu^3}\,,\nonumber\\
T&=&\displaystyle \pm\frac{2(\kappa')^2\tau+\kappa^2\tau^3+\kappa\kappa'\tau'-\kappa\kappa''\tau}{\nu^{\frac52}}\,. \nonumber
\end{eqnarray}
\end{prop}
\begin{demo}
We give a proof of the above equations in our context. 
We compute the right hand sides of \eqref{f_T2} and \eqref{f_Q_gamma_sigma} using the formulae \eqref{f_sigma} of $\sigma$, \eqref{f_gamma_XX} of $\gamma$, \eqref{f_drho_ds} of $\frac{d\rho}{ds}$, and \eqref{f_<>_vect_prod} and \eqref{f_<>_wedge_prod} of our pseudo-Riemannian structure . 
%
%
%
We use some data of $\langle d^{\,i}m/ds^i, d^{\,j}m/ds^j \rangle$, which are shown in table \ref{table_mij}, where $F_i=\big\|d^{\,i}m/ds^i\big\|^2$ are given by \setlength\arraycolsep{1pt}
$$
\begin{array}{rcl}
F_2&=&\kappa^2, \\
F_3&=&\kappa^4+{\kappa^{\prime}}^2+\kappa^2\tau^2, \\
F_4&=& 9\kappa^2 (\kappa')^2 + (\kappa^3 + \kappa \tau^2 -\kappa'')^2 +(2\kappa'\tau+\kappa\tau')^2. 
\end{array}
$$
We remark that $\nu=(\frac{d\rho}{ds})^2$ is given by $F_2$ and $F_3$ as $\nu^2 =F_3 -F_2^2$ (\cite{Li}). 

\begin{table}[h]
\begin{center}
\par\noindent
\begin{tabular}{|c|c|c|c|c|c|c|}
\hline
$\langle\cdot,\cdot\rangle$ $\phantom{\stackrel{{a}}{0}}$
& $m$  &  $m'$  &  $m''$  & $m'''$ & $m^{(4)}$ & $m^{(5)}$  \\[1mm] 
\hline
 $\!\!\phantom{\stackrel{{a}}{0}}\!\!m$ & $0$ & $0$ & $-1$ & $0$ & $F_2$ & $\frac{5}{2}F_2'$ \\[1mm] 
\hline
 $\!\!\phantom{\stackrel{{a}}{0}}\!\!m'$ &  & $1$ & $0$ & $-F_2$ & ${-\frac32{F_2}'}$ & $-2F_2''+F_3$  \\[1mm] 
\hline
 $\!\!\phantom{\stackrel{{a}}{0}}\!\!m''$ &  &  & $F_2$ & ${\frac12{F_2}'}$ & ${\frac12{F_2}''-F_3}$ & $\frac{1}{2}F_2'''-\frac{3}{2}F_3'$  \\[1mm] 
\hline
 $\!\!\phantom{\stackrel{{a}}{0}}\!\!m'''$ &  &  &  & $F_3$ & ${\frac12{F_3}'}$ & $\frac{1}{2}F_3''-F_4$  \\[1mm]
\hline 
 $\!\!\phantom{\stackrel{{a}}{0}}\!\!m^{(4)}$ & & & & & $F_4$ & $\frac{1}{2}F_4'$ \\[1mm]
\hline
\end{tabular}
\end{center}
\caption{A table of $\Big\langle \mbox{\large $\frac{d^{\,i}m}{{ds}^{i}}, \frac{d^{j}m}{{ds}^{j}} $}\Big\rangle$}
\label{table_mij}
\end{table}

Then we have 
\setlength\arraycolsep{1pt}
$$
\begin{array}{rcl}
%
\ccc\ga &=&\displaystyle \frac{2(\nu')^2-\nu\nu''}{2\nu^{\frac{7}{2}}}\ga' - \frac{3\nu'}{2\nu^{\frac{5}{2}}}\ga'' + \frac{1}{\nu^{\frac{3}{2}}}\ga''', \\[4mm]
\langle \ccc\ga, \ccc\ga \rangle &=&\displaystyle  \frac{(2(\nu')^2-\nu\nu'')^2}{4\nu^7} \langle \ga' , \ga' \rangle - \frac{3\nu'(2(\nu')^2-\nu\nu'')}{2\nu^6} \langle \ga' , \ga'' \rangle + \frac{2(\nu')^2-\nu\nu''}{\nu^5} \langle \ga' , \ga''' \rangle \\[4mm]
                                 &&\displaystyle +\frac{9(\nu')^2}{4\nu^5} \langle \ga'' , \ga'' \rangle - \frac{3\nu'}{\nu^4} \langle \ga'' , \ga''' \rangle +\frac{1}{\nu^3} \langle \ga''' , \ga''' \rangle \\[4mm]
                                 &=&\displaystyle \frac{-11(\nu')^2+4\nu\nu'' +4 \langle \ga''' , \ga''' \rangle }{4\nu^3},\\[4mm]
\langle \ga''', \ga''' \rangle &=& \displaystyle 5F_2^3 -11F_2F_3 +8F_2F_2'' - \frac{23}{2}F_2' -F_3'' +6F_4, \\[4mm]
\si' &=& \displaystyle \frac{1}{\nu} m \times m' \times m'' \times m^{(4)} - \frac{\nu'}{\nu^2} m \times m \times m' \times m'' \times m''', \\[4mm]
\langle \c\si , \c\si \rangle &=& \frac{1}{\nu} \langle \si', \si' \rangle \\[2mm]
                              &=&\displaystyle \frac{F_2^3 +F_2F_2''-2F_2F_3+F_4-\frac{9}{4}(F_2')^2-(\nu')^2}{\nu^3}\,.
\end{array}
$$
Now the conclusion follows from direct computation. 
\end{demo}

\section{Characterization of canal surfaces}
\subsection{Canal surfaces}
\begin{defi}\rm 
A {\em canal surface} in $\ss^3$ is the envelope of a space-like curve $\sigma$ in $\Lambda^4$. 
A circle $\spa(\sigma, \Pt{\sigma})^{\bot}$ is called a {\em characteristic circle} of the canal surface $\sigma$. 
The corresponding point in ${\mathcal S}(1,3)\subset \stackrel{2}{\bigwedge}\RR^5_1$ is given by $\gamma=\si\wedge\Pt\si$, where $\Pt{\phantom{\si}}$ denotes the derivative with respect to the arc-length parameter of $\sigma$. 
\end{defi}

There are three types of canal surfaces. 
Suppose a canal surface $\sigma$ is parameterized by the arc-length. 
The geodesic curvature vector $\vect{k_g}=\si+\Ppt\si$ can be one of the three types: time-like, space-like, or light-like. 
If $\vect{k_g}$ is always time-like, the canal surface is {\em regular}, that is an immersed cylinder. 
If $\vect{k_g}$ is always space-like, the canal surface has two singular curves. The characteristic circles are tangent to these two curves. 
If $\vect{k_g}$ is always light-like, the canal surface is formed by the osculating circles of a curve except for degenerates cases like families of spheres tangent to a given direction at a point. 

According to the three types of the canal surface, the corresponding curve of characteristic circles $\gamma$ is time-like, space-like, or light-like. 
This is shown by the following proposition. 

\begin{prop}
We have $\langle\vect{k_g}, \vect{k_g}\rangle=\langle\Pt\gamma,\Pt\gamma\rangle$. 
\end{prop}

\begin{demo}
Since $\langle\si,\si\rangle=1$ and $\langle\Pt\si, \Pt\si\rangle=1$, $\langle\si, \Pt\si\rangle=0$ and $\langle\si, \Ppt\si\rangle=-1$, which implies 
\[
\langle\vect{k_g}, \vect{k_g}\rangle
=\langle\si+\Ppt\si, \si+\Ppt\si\rangle
=\langle\Ppt\si, \Ppt\si\rangle-1
=\langle\si\wedge\Ppt\si, \si\wedge\Ppt\si\rangle
=\langle\Pt\gamma,\Pt\gamma\rangle\,.
\]
\end{demo}

The reader is referred to see \cite{La-So} for further study of canal surfaces. 

\subsection{Canal surfaces in terms of a curve in ${\mathcal S}(1,3)$}
\begin{theo}
A curve $\gamma$ in ${\mathcal S}(1,3)$ is a family of characteristic circles of a canal surface $\si$ if and only if the vector $\Pt{\gamma}$ always satisfies the Pl\"ucker relations  \eqref{Plucker_relation_1} - \eqref{Plucker_relation_5}. 
\end{theo}

\bigskip
\begin{demo}
(1) The ``only if" part. 
Suppose $\gamma$ is a curve of characteristic circles of a canal surface $\si$. 
Then $\gamma=\si\wedge\Pt\si$, where $\Pt{\phantom{\si}}$ denotes the derivative with respect to the arc-length parameter of $\sigma$. 
Therefore $\Pt\gamma=\sigma \wedge \Ppt{\sigma}$ is a pure vector. 

\smallskip
(2) The ``if'' part. 
Suppose a curve $\gamma$ can be expressed as $\gamma(t)=\alpha(t)\wedge\beta(t)$ for some $\alpha(t)$ and $\beta(t)$ in $\RR^5_1$ with $\langle\alpha(t),\alpha(t)\rangle=\esc{\beta(t)}{\beta(t)}=1$ and $\esc{\alpha(t)}{\beta(t)}=0$. 
The condition is equivalent to $\gamma'$ always satisfying the Pl\"ucker relation. 
We may assume without loss of generality that, at $t=t_0$, 
\[\alpha=\alpha(t_0)=(0,1,0,0,0) \>\>\mbox{ and } \>\>\beta=\beta(t_0)=(0,0,1,0,0).\]
Suppose $\alpha'=\alpha'(t_0)$ and $\beta'=\beta'(t_0)$ are given by 
\[\alpha'=(a,b,c,d,e)\>\>\mbox{ and } \>\>\beta'=(a',b',c',d',e').\]
Then $b=c'=0$ as $\esc{\alpha}{\alpha'}=\esc{\beta}{\beta'}=0$ and $b'+c=0$ as $\esc{\alpha}{\beta'}+\esc{\alpha'}{\beta}=0$. 
Using the formula \eqref{Plucker_coordinates}, 
the Pl\"ucker coordinates of $\gamma'=\alpha\wedge\beta'+\alpha'\wedge\beta$ are given by 
\[
(p_{34},\cdots,p_{12}\,;\,p_{04},\cdots,p_{01})(\alpha\wedge\beta'+\alpha'\wedge\beta)
=(\,0\,,-e,-d,e',d',\,0\,;\,0\,,\,0\,,-a,-a').
\]
As they must satisfy the Pl\"ucker relations (\ref{Plucker_relation_1} - \ref{Plucker_relation_5}) by the assumption, we have 
\[
a'd=ad', \, a'e=ae', \>\mbox{ and } \> d'e=de'.
\]
Put $u=(a',\,0\,,\,0\,,d',e')$ and $f={a}/{a'}={d}/{d'}={e}/{e'}$, 
then $\alpha'=fu+c\beta$ and $\beta'=u-c\alpha$, and therefore 
\[
\gamma'=\alpha\wedge\beta'+\alpha'\wedge\beta
=(\alpha-f\beta)\wedge u, 
\]
which implies that $\gamma$ and $\gamma'$ correspond to 
$\Pi=\spa(\alpha,\beta)$ and $P=\spa(\alpha-f\beta,u)$ respectively. 
It follows that $\Pi\cap P=\spa(\alpha-f\beta)$ is a space-like line. 

Let $\sigma$ be one of the unit vectors in $\Pi\cap P$: 
\[
\sigma=(\cos\theta)\alpha+(\sin\theta)\beta \hspace{1cm}(\theta=\arctan(-f)).
\]
Since 
\[
\sigma'=(\cos\theta)\alpha'+(\sin\theta)\beta'-\theta'(\sin\theta)\alpha+\theta'(\cos\theta)\beta
\]
and $\alpha'=fu+c\beta=-(\tan\theta)u+c\beta$  and $\beta'=u-c\alpha$ we have 
\[
\sigma\wedge\sigma'=(\theta'+c)\alpha\wedge\beta=(\theta'+c)\gamma,
\]
which completes the proof. 
\end{demo}

\begin{coro}
There is a bijection 
\[
\varphi:
\{\mbox{\rm curve }\si\subset\Lambda^4\,|\,\Pt\si\mbox{ \rm is space-like}\}\to\{\mbox{\rm curve }\gamma\subset\mathcal{S}(1,3)\,|\,\gamma'\mbox{ \rm is pure}\}
\]
which is given by $\varphi(\si)=\si\wedge\Pt\si$, where $\Pt{\phantom{\si}}$ denotes the derivative with respect to the arc-length parameter of $\sigma$. 
Its inverse is given by $\varphi^{-1}(\gamma)=\Lambda^4\cap\,\spa(\Pi\cap P)$, where $\Pi$ and $P$ are $2$ dimensional space-like vector subspaces of $\RR^5_1$ which correspond to $\gamma$ and $\gamma'$ respectively. 
\end{coro}

\bibliographystyle{plain}

\end{document}